# A positive solution to the Busemann-Petty problem in $\mathbb{R}^4$

By Gaoyong Zhang*

## Introduction

Motivated by basic questions in Minkowski geometry, H. Busemann and C. M. Petty posed ten problems about convex bodies in 1956 (see [BP]). The first problem, now known as the Busemann-Petty problem, states:

*If $K$ and $L$ are origin-symmetric convex bodies in $\mathbb{R}^n$, and for each hyperplane $H$ through the origin the volumes of their central slices satisfy*

$$\mathrm{vol}_{n-1}(K \cap H) < \mathrm{vol}_{n-1}(L \cap H),$$

*does it follow that the volumes of the bodies themselves satisfy*

$$\mathrm{vol}_n(K) < \mathrm{vol}_n(L)?$$

The problem is trivially positive in $\mathbb{R}^2$. However, a surprising negative answer for $n \geq 12$ was given by Larman and Rogers [LR] in 1975. Subsequently, a series of contributions were made to reduce the dimensions to $n \geq 5$ by a number of authors (see [Ba], [Bo], [G2], [Gi], [Pa], and [Z1]). That is, the problem has a negative answer for $n \geq 5$. See [G3] for a detailed description. It was proved by Gardner [G1] that the problem has a positive answer for $n = 3$. The case of $n = 4$ was considered in [Z1]. But the answer to this case in [Z1] is not correct. This paper presents the correct solution, namely, the Busemann-Petty problem has a positive solution in $\mathbb{R}^4$, which, together with results of other cases, brings the Busemann-Petty problem to a conclusion.

A key step to the solution of the Busemann-Petty problem is the discovery of the relation of the problem to intersection bodies by Lutwak [Lu]. An origin-symmetric convex body $K$ in $\mathbb{R}^n$ is called an intersection body if its radial function $\rho_K$ is the spherical Radon transform of a nonnegative measure $\mu$ on the unit sphere $S^{n-1}$. The value of the radial function of $K$, $\rho_K(u)$, in the direction $u \in S^{n-1}$, is defined as the distance from the center of $K$ to its boundary in that direction. When $\mu$ is a positive continuous function, $K$ is called the intersection body of a star body. The notion of intersection body was introduced by Lutwak [Lu] who proved that the Busemann-Petty problem has a positive answer if $K$ is an intersection body in $\mathbb{R}^n$. Based on this relation, a positive answer to the Busemann-Petty problem in $\mathbb{R}^3$ was given by Gardner [G1] who showed that all origin-symmetric convex bodies in $\mathbb{R}^3$ are intersection bodies.

---

*Research supported, in part, by NSF grant DMS-9803261.

The relation of the Busemann-Petty problem to intersection bodies proved by Lutwak can be formulated as: A negative answer to the Busemann-Petty problem is equivalent to the existence of convex nonintersection bodies (see [G2] and [Z2]). The author attempted in [Z1] to give a negative answer for all dimensions $\geq 4$ by trying to show that cubes in $\mathbb{R}^n$ ($n \geq 4$) are not intersection bodies (see Theorem 5.3 in [Z1]). However, there is an error in Lemma 5.1 of [Z1]. It affects only Theorems 5.3 and 5.4 there. The correct version of Theorem 5.3 is that no cube in $\mathbb{R}^n$ ($n > 4$) is an intersection body. This follows immediately from Theorem 6.1 of [Z1] which says that no generalized cylinder in $\mathbb{R}^n$ ($n > 4$) is an intersection body. Note that the proof of Theorem 6.1 in [Z1] holds for intersection bodies, although the definition of intersection body of a star body was the one used in [Z1]. Therefore, Theorem 5.4 in [Z1] should have stated: The Busemann-Petty problem has a negative solution in $\mathbb{R}^n$ for $n > 4$.

In his important work [K1], Koldobsky applied the Fourier transform to the study of intersection bodies. In [K2], he showed that cubes in $\mathbb{R}^4$ are intersection bodies. It was this result that exposed the error mentioned above and led to the present paper, which presents the correct solution to the Busemann-Petty problem in $\mathbb{R}^4$. One of the key ideas in the proof, previously employed by Gardner [G1], is the use of cylindrical coordinates in computing the inverse spherical Radon transform.

## 1. The inverse Radon transform on $S^3$ and intersection bodies in $\mathbb{R}^4$

The radial function $\rho_L$ of a star body $L$ is defined by

$$\rho_L(u) = \max\{r \geq 0 : ru \in L\}, \quad u \in S^{n-1}.$$

It is required in this paper that the radial function is continuous and even. For basic facts about star bodies and convex bodies, see [G3] and [S].

For a continuous function $f$ on $S^{n-1}$, the spherical Radon transform $\mathrm{R}f$ of $f$ is defined by

$$(\mathrm{R}f)(u) = \int_{S^{n-1} \cap u^\perp} f(v)dv, \quad u \in S^{n-1},$$

where $u^\perp$ is the $(n-1)$-dimensional subspace orthogonal to the unit vector $u$. Since the spherical Radon transform is self-adjoint, one can define the Radon transform $R\mu$ for a measure $\mu$ on $S^{n-1}$ by

$$\langle R\mu, f \rangle = \langle \mu, Rf \rangle.$$

The intersection body $IL$ of star body $L$ is defined by

$$\rho_{IL}(u) = \mathrm{vol}_{n-1}(L \cap u^{-1}) = R\left(\frac{1}{n-1}\rho_L^{n-1}\right)(u), \quad u \in S^{n-1}.$$

An origin-symmetric convex body $K$ is called the *intersection body of a star body* if there exists a star body $L$ so that $K = IL$. That is, the inverse spherical Radon transform $R^{-1}\rho_K$ is a positive continuous function. A slight extension of this definition is that an origin-symmetric convex body $K$ is called an *intersection body* if the inverse spherical Radon transform $R^{-1}\rho_K$ is a non-negative measure.

Let $\Delta$ be the Laplacian on the unit sphere $S^3$. Helgason's inversion formula for the Radon transform $R$ on $S^3$ is (see [H, p. 161])

$$\frac{1}{16\pi^2}(1-\Delta)RR = 1.$$

It implies that

$$(1) \qquad R^{-1}\rho_K = \frac{1}{16\pi^2}R(1-\Delta)\rho_K$$

for an origin-symmetric convex body $K$ in $\mathbb{R}^4$. This formula shows that $R^{-1}\rho_K$ is continuous when $\rho_K$ is of class $C^2$. The following lemma provides an inversion formula which gives the positivity of $R^{-1}\rho_K$.

Let $K$ be an origin-symmetric convex body in $\mathbb{R}^4$, and let $A_u(z)$ be the volume of $K \cap (zu + u^\perp)$, where $z$ is real and $u \in S^3$.

LEMMA 1. *If $K$ is an origin-symmetric convex body in $\mathbb{R}^4$ whose boundary is of class $C^2$, then*

$$(2) \qquad (R^{-1}\rho_K)(u) = -\frac{1}{16\pi^2}A_u''(0), \quad u \in S^3.$$

*Proof.* By rotation, it suffices to prove (2) for the north pole of $S^3$. From Helgason's inversion formula (1), the inverse spherical Radon transform of $\rho_K$, $f = R^{-1}\rho_K$, is a continuous function when $\rho_K$ is of class $C^2$. Let

$$u = u(v, \phi) = (v\sin\phi, \cos\phi), \quad u \in S^3, \quad v \in S^2, \quad 0 \le \phi \le \pi,$$

and let $\rho_K(v, \phi) = \rho_K(u)$ be the radial function of $K$. Define

$$\bar{\rho}_K(\phi) = \int_{S^2} \rho_K(v, \phi) dv,$$

$$\bar{f}(\phi) = \int_{S^2} f(u) dv,$$

$$r(v, \phi) = \rho_K(v, \phi) \sin \phi,$$

$$\bar{r}(\phi) = \bar{\rho}_K(\phi) \sin \phi.$$

Consider $\bar{\rho}_K$ and $\bar{f}$ as functions on $S^3$ which are SO(3) invariant. Since the spherical Radon transform is intertwining, we have $\bar{\rho}_K = R\bar{f}$ (for a simple proof, see [G3, Th C.2.8]). From this and Lemma 2.1 in [Z1], or Theorem C.2.9 in [G3], we obtain

$$\bar{\rho}_K(\phi) = \frac{4\pi}{\sin \phi} \int_{\frac{\pi}{2}-\phi}^{\frac{\pi}{2}} \bar{f}(\psi) \sin \psi d\psi.$$

Taking the derivative on both sides of this equation gives

$$(\bar{\rho}_K(\phi) \sin \phi)' = 4\pi \bar{f}(\frac{\pi}{2} - \phi) \sin(\frac{\pi}{2} - \phi).$$

It follows that

$$4\pi \bar{f}(0) = \lim_{\phi \to \frac{\pi}{2}} \frac{(\bar{\rho}_K(\phi) \sin \phi)'}{\cos \phi} = -\bar{r}''(\frac{\pi}{2}).$$

Since $\frac{1}{4\pi} \bar{f}(0)$ is the value of $f$ at the north pole, we obtain

$$(3) \qquad\qquad f(u_0) = -\frac{1}{16\pi^2} \bar{r}''(\frac{\pi}{2}),$$

where $u_0$ is the north pole of $S^3$.

Consider the variable $z$ defined by $z = \rho_K \cos \phi$. Then $\tan \phi = \frac{r}{z}$. Differentiating this equation and using $\frac{1}{\cos^2 \phi} = 1 + \tan^2 \phi = 1 + \frac{r^2}{z^2}$ give

$$(4) \qquad\qquad z^2 + r^2 = z \frac{dr}{d\phi} - r \frac{dz}{d\phi}.$$

This yields

$$(5) \qquad\qquad \frac{dz}{d\phi}\Big|_{\phi=\frac{\pi}{2}} = -r(v, \frac{\pi}{2}).$$

Differentiating (4) gives

$$(6) \qquad\qquad 2z \frac{dz}{d\phi} + 2r \frac{dr}{d\phi} = z \frac{d^2r}{d\phi^2} - r \frac{d^2z}{d\phi^2}.$$

From (5),

$$(7) \qquad \frac{dr}{d\phi}\bigg|_{\phi=\frac{\pi}{2}} = \frac{dr}{dz}\frac{dz}{d\phi}\bigg|_{\phi=\frac{\pi}{2}} = -r\,\frac{dr}{dz}\bigg|_{z=0}.$$

From (6) and (7),

$$(8) \qquad \frac{d^2z}{d\phi^2}\bigg|_{\phi=\frac{\pi}{2}} = 2r\,\frac{dr}{dz}\bigg|_{z=0}.$$

From (5), (8), and

$$\frac{d^2r}{d\phi^2} = \frac{d^2r}{dz^2}\left(\frac{dz}{d\phi}\right)^2 + \frac{dr}{dz}\frac{d^2z}{d\phi^2},$$

we have

$$(9) \qquad \frac{d^2r}{d\phi^2}\bigg|_{\phi=\frac{\pi}{2}} = \frac{d^2r}{dz^2}\bigg|_{z=0} r(v,\tfrac{\pi}{2})^2 + 2r(v,\tfrac{\pi}{2})\left(\frac{dr}{dz}\right)^2_{z=0}$$

$$= \left(r^2\frac{d^2r}{dz^2}\right)_{z=0} + \left(2r\left(\frac{dr}{dz}\right)^2\right)_{z=0}$$

$$= \frac{1}{3}\,\frac{d^2r^3}{dz^2}\bigg|_{z=0}.$$

Integrating both sides of (9) over $S^2$ with respect to $v$ gives

$$\int_{S^2} \frac{d^2r}{d\phi^2}(v,\phi)\bigg|_{\phi=\frac{\pi}{2}} dv = \frac{1}{3}\int_{S^2}\frac{d^2r^3}{dz^2}(v,z)\bigg|_{z=0} dv.$$

Since $K$ has $C^2$ boundary, one can interchange the second order derivative and the integral on each side of the last equation. We obtain

$$\frac{d^2}{d\phi^2}\bar{r}(\phi)\bigg|_{\phi=\frac{\pi}{2}} = \frac{d^2}{dz^2}\left(\frac{1}{3}\int_{S^2}r^3(v,z)dv\right)_{z=0}.$$

Note that the 3-dimensional volume of the intersection of the hyperplane $x_4 = z$ with the convex body $K$, denoted by $A_{u_0}(z)$, is given by

$$A_{u_0}(z) = \frac{1}{3}\int_{S^2}r^3(v,z)dv.$$

Therefore, we have

$$(10) \qquad \bar{r}''(\tfrac{\pi}{2}) = A_{u_0}''(0).$$

Formula (2) follows from (3) and (10). $\qquad\qquad\qquad\qquad\qquad\square$

Recently, Gardner, Koldobsky and Schlumprecht [GKS] have generalized the formula (2) to $n$ dimensions by using techniques of the Fourier transform. A different proof of their formulas is given by Barthe, Fradelizi and Maurey [BFM].

THEOREM 2. *If $K$ is an origin-symmetric convex body in $\mathbb{R}^4$ whose boundary is of class $C^2$ and has positive curvature, then $K$ is an intersection body of a star body.*

*Proof.* By the Brunn-Minkowski inequality and the strict convexity of $K$, $A(z)^{\frac{1}{3}}$ is strictly concave. When one slices a symmetric convex body by parallel hyperplanes, the central section has maximal volume. Hence, $A'(0) = 0$. It follows that

$$A''(0) = 3A(0)^{\frac{2}{3}}\left(A(z)^{\frac{1}{3}}\right)''_{z=0} < 0.$$

By Lemma 1, $\mathrm{R}^{-1}\rho_K$ is a positive continuous function. Therefore, $K$ is the intersection body of a star body.                          □

When a convex body is identified with its radial function, the class of intersection bodies is closed under the uniform topology. Since every origin-symmetric convex body can be approximated by origin-symmetric convex bodies whose boundaries are of class $C^2$ and have positive curvatures, we obtain:

THEOREM 3. *All origin-symmetric convex bodies in $\mathbb{R}^4$ are intersection bodies.*

Theorem 3 is proved for convex bodies of revolution by Gardner [G2] and by Zhang [Z1], and is proved for cubes and other special cases by Koldobsky [K2]. In higher dimensions, the situation is different. For example, it is proved by Zhang [Z1] that generalized cylinders in $\mathbb{R}^n$ ($n > 4$) are not intersection bodies, and is proved by Koldobsky [K1] that the unit balls of finite dimensional subspaces of an $L_p$ space, $1 \leq p \leq 2$, are intersection bodies. In three dimensions, Gardner [G1] proved that all origin-symmetric convex bodies in $\mathbb{R}^3$ are intersection bodies. One can also prove this by Theorem 3 and a result of Fallert, Goodey and Weil [FGW] which says that central sections of intersection bodies are again intersection bodies. An intersection body may not be the intersection body of a star body. It is shown by Zhang [Z4] that no polytope in $\mathbb{R}^n$ ($n > 3$) is an intersection body of a star body. Campi [C] is able to prove a complete result which says that no polytope in $\mathbb{R}^n$ ($n > 2$) is an intersection body of a star body.

## 2. A positive solution to the Busemann-Petty problem in $\mathbb{R}^4$

The following relation of the Busemann-Petty problem to intersection bodies was proved by Lutwak [Lu].

THEOREM 4 (Lutwak). *The Busemann-Petty problem has a positive solution if the convex body with smaller cross sections is an intersection body.*

From Theorems 3 and 4, we conclude:

THEOREM 5. *The Busemann-Petty problem in $\mathbb{R}^4$ has a positive solution.*

From Theorem 3 and Corollary 2.19 in [Z2], we have the following corollary about the maximal cross section of a convex body.

COROLLARY 6. *If $K$ is an origin-symmetric convex body in $\mathbb{R}^4$, then*

$$(11) \qquad \mathrm{vol}_4(K)^{\frac{3}{4}} \leq \frac{3}{8}(\sqrt{2}\pi)^{\frac{1}{2}} \max_{u \in S^3} \mathrm{vol}_3(K \cap u^\perp)$$

*with equality if and only if $K$ is a ball.*

Inequality (11) implies that, in $\mathbb{R}^4$, balls attain the minmax of the volume of central hyperplane sections of origin-symmetric convex bodies with fixed volume. The corresponding inequality in $\mathbb{R}^3$ to inequality (11) was proved by Gardner (see [G3, Th. 9.4.11]). However, it is no longer the case in higher dimensions at least for $n \geq 7$. Ball [Ba] showed that cubes are counterexamples for $n \geq 10$. Giannopoulos [Gi] showed that certain cylinders are counterexamples for $n \geq 7$. The following question, known as the slicing problem, has been of interest (see [MP] for details):

*Does there exist a positive constant $c$ independent of the dimension $n$ so that*

$$\mathrm{vol}_n(K)^{\frac{n-1}{n}} \leq c \max_{u \in S^{n-1}} \mathrm{vol}_{n-1}(K \cap u^\perp)$$

*for every origin-symmetric convex body $K$ in $\mathbb{R}^n$?*

## 3. The generalized Busemann-Petty problem

Besides considering hyperplane sections, one can also consider intermediate sections of convex bodies. For a fixed integer $1 < i < n$, the Busemann-Petty problem has the following generalization (see Problem 8.2 in [G3]):

*If $K$ and $L$ are origin-symmetric convex bodies in $\mathbb{R}^n$, and for every $i$-dimensional subspace $H$ the volumes of sections satisfy*

$$\mathrm{vol}_i(K \cap H) < \mathrm{vol}_i(L \cap H),$$

*does it follow that the volumes of the bodies themselves satisfy*

$$\mathrm{vol}_n(K) < \mathrm{vol}_n(L)?$$

When $i = n-1$, this is the Busemann-Petty problem. It turns out that the solution to the generalized Busemann-Petty problem depends strongly on the dimension $i$ of the sections of convex bodies. It is proved by Bourgain and Zhang [BoZ] that the solution is negative when $3 < i < n$. The generalized Busemann-Petty problem has a positive solution when $K$ belongs to a certain class of convex bodies, called $i$-intersection bodies, which contains all intersection bodies (see Theorem 5 in [Z3] and Lemma 6.1 in [GrZ]). In particular, when $K$ is an intersection body, the generalized Busemann-Petty problem has a positive solution. From this fact and Theorem 3, we have:

THEOREM 7. *The generalized Busemann-Petty problem in $\mathbb{R}^4$ has a positive solution.*

It might be still true that the generalized Busemann-Petty problem has a positive solution when $i = 2, 3$, and $n \geq 5$. This remains open.

*Acknowledgement.* I am very grateful to Professors R. J. Gardner, E. Grinberg, and E. Lutwak for their encouragement while this work was done.

POLYTECHNIC UNIVERSITY, BROOKLYN, NY
*E-mail address:* gzhang@math.poly.edu